\documentclass[10pt]{amsart}
\newcommand{\vs}[1]{\vspace{#1}}

\newcommand{\ra}{\rightarrow}
\newcommand{\R}{\mathbb{R}}
\newcommand{\T}{\mathbb{T}}
\newcommand{\Z}{\mathbb{Z}}
\newcommand{\N}{\mathbb{N}}

\newcommand{\ep}{\epsilon}
\newcommand{\om}{\omega}

\newcommand{\de}{\delta}

\newcommand{\8}{\infty}
\newcommand{\nn}{\nonumber}
\newcommand{\be}{\begin{eqnarray}}
\newcommand{\ee}{\end{eqnarray}}
\newcommand{\dfr}{\mbox{{\rm Diff}$_\omega^r(M)$}}
\newcommand{\dfro}{\mbox{{\rm Diff}$_{\omega, 0}^r(M)$}}

\newcommand{\bv}{\mbox{{\bf v}}}
\newtheorem{thm}{Theorem}[section]
\newtheorem{lem}[thm]{Lemma}

\newtheorem{cor}[thm]{Corollary}
\usepackage{graphicx}
\begin{document}

\title[Homoclinic point]{Homoclinic Points For Area-Preserving Surface
  Diffeomorphisms}
\author{Zhihong Xia}
\address{Department of Mathematics \\ Northwestern University \\
Evanston, Illinois 60208}
\thanks{Research supported in part by National Science Foundation.}
\date{June 10, 2006, daft version}
\email{xia@math.northwestern.edu}

\begin{abstract}
  We show a $C^r$ connecting lemma for area-preserving surface
  diffeomorphisms and for periodic Hamiltonian on surfaces. We prove
  that for a generic $C^r$, $r=1, 2, \ldots$, $\8$, area-preserving
  diffeomorphism on a compact orientable surface, homotopic to
  identity, every hyperbolic periodic point has a transversal
  homoclinic point. We also show that for a $C^r$, $r=1, 2, \ldots$,
  $\8$ generic time periodic Hamiltonian vector field in a compact
  orientable surface, every hyperbolic periodic trajectory has a
  transversal homoclinic point. The proof explores the special
  properties of diffeomorphisms that are generated by Hamiltonian
  flows.
\end{abstract}

\maketitle

\section{Introduction and statement of main results}

Let $M^{2n}$ be a compact $2n$ dimensional symplectic manifold with a
symplectic form $\om$ and let $\dfr$ be the set of all $C^r$, $r=1, 2,
\cdots, \8$, diffeomorphisms that preserves the symplectic form
$\om$. For $f \in \dfr$, a point $p \in M$ is said to be a periodic
point of $f$ with period $k$ if $f^k(p)=p$. A periodic point $p$ is
said to be hyperbolic if all the eigenvalues of $df^k(p)$ are away
from the unit circle. For any hyperbolic periodic point, there stable
manifold, where points approaches $p$ under forward iterations of
$f^k$, and a unstable manifold, where points approaches $p$ under
backward iterations of $f^k$. For symplectic diffeomorphisms, both the
stable and unstable manifolds are $n$ dimensional. We write the stable
and unstable manifold of $p$ as, respectively, $W^s_f(p)$ and
$W^u_f(p)$.

A point $q \in M$ is said to be a homoclinic point to a hyperbolic
periodic point $p$ if $q \in (W^s(p) \cap W^u(p)) \backslash
\{p\}$. i.e., homoclinic points are nontrivial intersections of stable
and unstable manifolds. A homoclinic point $q$ is said to be
transversal if $W^s_f(p)$ intersects $W^u_f(p)$ transversally at $q$.

Transversal homoclinic points are responsible for very complicated and
chaotic dynamics. Poincar\'e discovered the homoclinic phenomenon and
its associated chaotic dynamics in his study of the restricted three
body problem in celestial mechanics \cite{Po1892}. Poincar\'e
conjectured that transversal homoclinic points occurs generically in
Hamiltonian systems. Proving the generic existence of homoclinic
points to a hyperbolic periodic trajectory is a classical problem in
Hamiltonian dynamics.  Our main result in this paper positively
answers Poincar\'e's question for lower dimensional Hamiltonian
systems.

Our first result is on area-preserving diffeomorphisms on compact
surfaces that are homotopic to identity. Recall that a subset is said
to be {\it residual}\/ if it contains a countable intersection of open
and dense subset. A property is said to be {\it generic}\/ if it holds
on a residual set.

\begin{thm}
  Let $M$ be a compact orientable surface and let $\dfro$ be the set of
  all $C^r$ diffeomorphisms on $M$ preserving an area form $\omega$ on
  $M$ and isotropic to identity. Then for any $r =1,2, \ldots, \8$,
  there is a residual subset $R_1 \subset \dfro$ such that if $f \in
  R_1$ and $p$ is a hyperbolic periodic point, then $$W^s(p) \cap
  W^u(p) \backslash \{p\} \neq \emptyset$$
\label{mainthm}\end{thm}

Our result naturally applies to Hamiltonian systems. Again let
$M^{2n}$ be a symplectic manifold with the symplectic form $\om$. The
nondegenerate closed two-form $\om$ on $M^{2n}$ defines an isomorphism
between the tangent bundle and the cotangent bundle of $M$ with $J:
T^*M \ra TM$, where the isomorphism $J$ is uniquely determined by the
following equality: $\om(*, J\alpha) = \alpha(*)$ for any one-form
$\alpha$.

Let $H: M^{2n} \ra \R$ be a $C^{r+1}$ real valued function. $JdH$
defines a $C^r$ vector field on $M$. This vector field is call a
Hamiltonian vector field and the function $H$ is called a Hamiltonian
function, or Hamiltonian. The Hamiltonian is a constant of motion
under the Hamiltonian flow. A periodic trajectory, also called closed
trajectory is said to be hyperbolic is all its characteristic
multipliers are away from the unit circle, except two. The
characteristic multipliers are one in the flow direction and in the
normal direction of the Hamiltonian function. The stable manifold and
unstable manifold of a hyperbolic closed trajectory both have
dimension $n$.

We will consider time dependent Hamiltonian systems, where the
Hamiltonian function can be written as $H = H(x, t)$, with $x \in M$
and $t \in \R$. The resulting Hamiltonian vector field is no
longer autonomous. We are particularly interested in time periodic
Hamiltonian systems where there is a positive real number $T$ such
that $H(x, t+T) \equiv H(x, t)$ for all $x \in M$ and $t \in \R$. Time
periodic Hamiltonians define Hamiltonian flow on $M \times S^1$.
Time periodic Hamiltonian can be reduced to a symplectic
diffeomorphisms on $M$ by considering its Poincar\'e map. Hyperbolic
periodic points and their stable and unstable manifold can be defined
in the same way as those of symplectic diffeomorphisms.

We can now state our theorem for the Hamiltonian case.

\begin{thm}
  Let $M$ be a compact orientable surface with an area form
  $\omega$. For any $r =1,2, \ldots, \8$, there is a residual subset
  $R_2 \subset C^{r+1}(M \times S^1)$ such that if $H \in R_2$ and
  $\gamma$ is a hyperbolic periodic orbit for the time periodic
  Hamiltonian vector field $JdH(x, t)$, then $$W^s(\gamma) \cap
  W^u(\gamma) \backslash \{ \gamma\} \neq \emptyset.$$
\label{mainthmh}\end{thm}

Related to the existence of homoclinic points in the Hamiltonian
systems and symplectic diffeomorphisms is the so called the closing
lemma, also raised by Poincar\'e \cite{Po1892} in the same
time. Poincar\'e believed that the periodic points are dense in a
typical Hamiltonian system on compact manifolds. For $C^1$ generic
diffeomorphisms ($C^2$ generic Hamiltonian systems), the closing lemma
was proved by Pugh \cite{Pu67} and Pugh and Robinson \cite{PR83}. A
different proof was given by Liao \cite{Liao79}]. The proof uses local
perturbation methods.  However, the $C^r$ closing lemma for arbitrary
$r$ remains open except for some very special cases (Anosov
\cite{Anosov67} for Anosov diffeomorphisms, Pesin \cite{Pesin77} for
non-uniformly hyperbolic diffeomorphisms, Xia and Zhang \cite{XZ05}
for some special cases of partially hyperbolic diffeomorphisms). Franks
and Le Calvez \cite{FL03}, Xia \cite{Xia04b} provided some further
evidence supporting $C^r$ closing lemma for area-preserving surface
diffeomorphisms.

Proving the existence of homoclinic points, or connecting stable and
unstable manifolds, is also know as a connecting lemma. Takens
\cite{Takens72} proved that $C^1$ generically every hyperbolic
periodic points has a transversal homoclinic point for symplectic and
volume-preserving diffeomorphisms. Takens results were extended by Xia
\cite{Xi96a}, using closing lemma types of techniques of
Hayashi. Hayashi established the first connecting lemma for hyperbolic
invariant set of general diffeomorphisms (Hayashi \cite{Ha97}). A
general version of $C^1$ connecting lemma was obtained by Wen and Xia
\cite{WX99} \cite{WX00}. We remark that he method for all these $C^1$
results are local perturbation methods and they are definitely
restricted to $C^1$ topology (cf. Gutierrrez \cite{Gutierrez87})

The $C^r$ connecting lemma for $r >1$ is a much more difficult
problem. Using the idea of ``closing gates'', Robinson
\cite{Robinson73}, following an idea of Newhouse, showed that one can
connect the stable and the unstable manifolds of a hyperbolic fixed
point on two sphere by a $C^r$ small perturbation, if the stable
manifold accumulates on the unstable manifold. The accumulation
condition is a generic condition for area-preserving
diffeomorphisms. Pixton \cite{Pixton82} extended Robinson's result to
hyperbolic periodic points on two sphere. Using a more topological
approach, Oliveira \cite{Oliveira87} showed the $C^r$ generic
existence of homoclinic points for area-preserving diffeomorphisms on
two torus $T^2$. More recently, Oliveira \cite{Oliveira00} showed the
generic existence of homoclinic points on surfaces of higher genus for
certain homotopic classes of area-preserving diffeomorphisms, where
the actions on first homology, if complicated enough, simply forces
intersections of stable and unstable manifolds. His result does not
apply to the important cases, including Hamiltonians, where the
action on homology is trivial.

To obtain our results, we explore the special properties of the
Hamiltonian flow and Hamiltonian diffeomorphisms. Our result uses the
concept of flux, a dual concept to the mean rotation vectors, for area
preserving diffeomorphisms that are homotopic to identity. We show
that the maps with rational flux has certain special properties. These
special properties enable us to show the existence of homoclinic
points.

We are motivated by the Arnold conjecture for Hamiltonian
diffeomorphisms on compact surfaces. Hamiltonian diffeomorphisms are
the ones generated by time periodic Hamiltonian flow. The number of
fixed points for these diffeomorphisms are larger than what is
predicted by the Lefschitz fixed point theorem, as Arnold conjectured
and later proved in this case by Floer \cite{Floer86}. Likewise, we
show that typically the stable and unstable manifolds for hyperbolic
fixed point do intersect. We remark that there are easy examples of
non-Hamiltonian diffeomorphisms with stable manifold and unstable
manifold accumulating on each other, but never intersecting.

Combining with generic existence of hyperbolic periodic points for
area-preserving surface diffeomorphisms (cf.\ Xia \cite{Xia04b}), we
have the following corollary. 

\begin{cor}
For $r=1, 2, \ldots, \8$, an open and dense set of $C^{r+1}$ time
periodic Hamiltonian systems on compact surfaces have positive
topological entropy.
\end{cor}

Similar statement is true for area-preserving diffeomorphisms in
$\dfro$.

\section{Flux and mean rotation numbers}
Let $M_g$, $g \geq 1$, be the orientable compact surface of genus $g$
and let $(a_i, b_i)$, $i=1, 2, \ldots, k$, be the canonical generators
of its first homology $H_1(M_g, \R)$ and let $(\alpha_i, \beta_j)$,
$i=1, 2, \ldots, k$, be the dual basis for the first De
Rham-cohomology $H^1(M_g, \R)$. By normalize the area form $\omega$ on
$M_g$, we may assume that $$\int_M \omega =1.$$

We consider the cases where $f$ is homotopic to the identity map. Let
$l$ be an oriented closed curve in $M = M_g$, since $f$ is isotopic to
identity, there is a oriented disk $D \subset M$ such that the
boundary $\partial D = f(l) - l$. We define the flux of $f$ across $l$
to be $$F_f(l) = \int_D d\mu \mod 1,$$ Where $\mu$ is the area element
given by the two form $\omega$. The above quantity is independent of
the choice of $D$. If $D'$ is another disk such that $\partial D' =
f(l)-l$, then $\partial (D-D') = \partial D - \partial D' = 0$. Since
$M$ is two dimensional, the boundaryless disk $D-D'$ is either the
whole manifold or the empty set. Therefore $\int_D d\mu = \int_{D'}
d\mu \mod 1$. We remark that one can always choose $D'=-(M\backslash
D)$, therefore the flux can only be defined up to mod 1.

In fact, the flux $F_f(l)$ depends only on the homology class of
$l$. i.e., if $l'$ is homologous to $l$, then $F_f(l') = F_f(l)$. This is
because that there will be is a disk $A \subset M$ such that $\partial
A = l-l'$. Let $D$ and $D'$ be the disks such that $\partial D = f(l)
-l$ and $\partial D' = f(l') -l'$, then $\partial f(A) = f (\partial
A) = f(l) - f(l')$ and $\partial D - \partial D' = \partial f(A) -
\partial A$, therefore, $D-D' = f(A) -A$. Consequently, $\int_D d\mu
-\int_{D'} d\mu = \int_{f(A)} d\mu - \int_A d\mu=0$.

The flux $F_f$ defined above can be extended to a linear map on the
first homology
$$F_f: H_1(M) \ra \R \mod 1,$$
and thus it can be represented by a cohomology vector.

The flux has a nice additive property under the iterations of the map
$f$. One easily verifies that $F_{f^k}(l) = kF_f(l), \mod 1$, for all
integers $k$ and for all closed curves $l$.

The flux is closely related to the {\it mean rotation vector}\/ of
$f$. Let $f \in \dfro$ be a area-preserving map homotopic to identity
and let $F_s: M \ra M$, $s \in [0, 1]$ be the homotopy: $F_0 = Id_M$,
the identity map on $M$, and $F_1 =f$. Let $\tilde{F}_s$ be a lift of
$F_s$ in the universal covering space $\tilde{M}$ of $M$. Let $\alpha$
be a differential one-form on $M$. To simplify the notation, we also
use $\alpha$ to denote its pull-back on $\tilde{M}$.  Let
$\tilde{M}_0$ be a fundamental domain in $\tilde{M}$ and let
$$R(\alpha) = \int_{\tilde{M}_0} (\int_{\{\tilde{F}_s(x): s \in [0,
  1]\} } \alpha) d\mu.$$ The function $R(\alpha)$ is linear on
$\alpha$. One can easily show that, when $\alpha$ is exact,
$R(\alpha)$ is zero. Therefore $R(\alpha)$ induces a linear map on the
first cohomology of $M$, $H^1(M, \R)$. This gives a homology vector in
$H_1(M, \R)$. This homology vector is called the mean rotation
vector. A different lift of $F_s$ will give a different mean rotation
vector $R(\alpha)$, but the difference is an integer.

We remark that the flux is more intuitive than the mean rotation
vector and it can be easily extended to the some cases where the map
is not isotopic to identity. For example, if a closed curve is
homologous to its image, then then flux across this closed curve can be
defined in the same way.

Let $\bv_f = (F_f(a_1), F_f(b_1), F_f(a_2), F_f(b_2), \ldots,
F_f(a_g), F_f(b_g))$ be the vector in $\T^{2g}$. We call $\bv_f$ the
flux vector. 

Let $f \in \dfr$ be a map with a rational flux vector, where all its
components are rational numbers. Then there exists
an integer $k$ such that $f^k$ has zero as its flux vector. We will
show that any map can be approximated by a map with rational flux
vector. Later in the paper, we will discuss special properties of maps
with zero flux vector.

It is a easy exercise to show that if a map is isotopic to identity
and has flux zero, then its mean rotation vector is zero, modulus
$\Z^n$.

As we will see later, all maps defined by Hamiltonian flows have zero
flux.

\section{Basic perturbations and some generic properties}

In this section, we will do a sequence of initial perturbations so
that the diffeomorphisms we consider satisfy certain properties. These
properties will later guarantee the existence of homoclinic points for
every hyperbolic periodic point.

Our first perturbation is to change the flux vector of the maps that
are homotopic to identity. We will show that the maps with rational
flux vector form a dense subset in $\dfr$.

\begin{lem}\label{lem70}
Let $f \in \dfr$ be a area preserving diffeomorphism, homotopic to
identity. Then for any neighborhood $V$ of $f$ in $\dfr$, there is a
map $g$ in $V$ such that the flux vector $\bv_g$ is rational.
\end{lem}

\begin{proof}
  Let $(a_i, b_i)$ be the canonical generators of the first homology
  $H_1(M, \Z)$. For any $i=1, 2, \ldots, g$, we may assume that $a_i$
  and $b_i$ are simple closed curves such that these curves don't
  intersect each other except $a_i$ and $b_i$; $a_i$ and $b_i$
  intersect at only one point and the intersection is transversal.

  Fix $i$, let $\delta_{b_i}$ be a small tubular neighborhood of
  $b_i$. We can parametrize this tubular neighborhood $\delta_{b_i}$
  by $\delta_{b_i}: S^1 \times [-\delta, \delta] \ra M$ for some small
  $\delta>0$. In fact, for convenience, we can even assume, without
  loss of generality, that the parametrization $\delta_{b_i}$ is
  area-preserving.

  Let $\beta: [-\delta, \delta] \rightarrow \mathbb{R}$ be a
  $C^\infty$ function such that $\beta(t) >0$ for all $-\delta < t <
  \delta$ and $\beta(-\delta) = \beta(\delta) =0$ and $\beta$ is
  $C^\infty$ flat at $\pm \delta$. i.e., all the derivatives of
  $\beta(t)$ at $\pm \delta$ are zero.

  Let $h_\ep: M \rightarrow AM$ be a $C^\infty$ diffeomorphism such
  that if $z \notin \delta_{b_i}$, $h_\ep(z) =z$ and if $z \in
  \delta_{b_i}$, $h_\ep = \delta_{a_i} \circ T_\ep \circ
  (\delta_{a_i})^{-1}$ where $T_\ep ( \theta, t) = (\theta + \ep
  \beta(t), t)$ for all $\theta \in S^1$ and $t \in [-\delta,
  \delta]$. We remark that $h_\ep \rightarrow \mbox{Id}_M$ in
  $C^\infty$ topology as $\ep \rightarrow 0$.

It is easy to see that the flux of $h_\ep$ across $a_j$ and $b_j$  are
all zero for $j \neq i$ and the flux across $b_i$ is also zero. For
$\ep>0$ small enough, we have $F_{h_\ep}(a_i) \neq 0$.

Since $F_{f \circ h_\ep}(l) = F_{f}(l) + F_{h_\ep}(l)$ for any closed
curve $l$, we have that for sufficiently small $\ep >0$, $F_{f \circ
  h_\ep}(a_i) - F_f(a_i) \neq 0$. There are infinitely many choices of
small $\ep$ such that $F_{f \circ h_\ep}(a_i)$ is rational.

We can do similar perturbations to every closed curves $a_i$ and $b_i$
so that the flux across $b_i$ and $a_i$ are all rational.

This proves the lemma.
\end{proof}

We now state a well-known local perturbation lemma (cf. Newhouse
\cite{Ne77}, Robinson \cite{Robinson73}, Takens \cite{Takens72} and
Xia \cite{Xi96a}).

Let $d$ be a metric on $M^n$ induced from some Riemann structure and let
$B_\delta(x)$ denote the set of $y\in M^n$ with $d(x, y) < \delta$. We also
let $\bar{B}_\delta(x)$ denote the closure of $B_\delta(x)$.

\begin{lem} \label{lem3}
  (perturbation lemma) Let $M^n$ be an $n$-dimensional compact manifold.
 Fix $\phi \in \dfr$, $r \geq 1$. There exist constants $\ep_0> 0$ and
 $c>0$, depending on $\phi$, such that for any $x\in M^n$, and any $\psi
 \in \dfr$ such that $\| \phi - \psi\|_{C^r} < \ep_0$, and any positive
 numbers $0< \de \leq \ep_0$, $0 < \ep \leq \ep_0$, the following facts
 hold.
\begin{quote}
  if $d(y, x) < c \de^r\ep$, then there is a $\psi_1 \in \dfr$,
  $\| \psi_1 -\psi\|_{C^r} < \ep$ such that $\psi_1 \psi^{-1}(x) =y$,
  $\psi_1(z) = \psi(z)$ for all $z \notin \psi^{-1}(B_\de(x))$, and
  $\psi^{-1}(z) = \psi^{-1}_1(z)$ for all $z \notin B_\de(x)$.
\end{quote}
\end{lem}

The proof of this lemma uses the generating functions, which provide
a convenient tool in the study of Hamiltonian systems and symplectic
diffeomorphisms.

We remark that the local perturbation provided in the perturbation
lemma does not change the flux vector or the mean rotation vector of
the map.

Now we use the perturbation lemma to prove the following simple
result:

\begin{lem} \label{lem2}
  Let $\phi \in \dfr$ and let $p\in M^n$ be a hyperbolic periodic points of
  period $k$, with respect to $\phi$. For any $\ep >0$, any $q\in
  W^u_\phi(p)$ and any neighborhood $U$ of $q$, there exist a $\phi' \in
  \dfr$, $\|\phi - \phi'\|_{C^r} < \ep$ such that
  \begin{enumerate}
  \item Support$(\phi-\phi')\subset U$ and hence, $p$ is a hyperbolic
    periodic point of period $k$ for $\phi'$,
  \item $q \in W^u_{\phi'}(p)$,
  \item $q$ is a recurrent point under $\phi'$.
  \end{enumerate}
\end{lem}

Recall that a point $q \in M^n$ is a {\em recurrent} point under $\phi'$ if
there exists a sequence of positive integers $\{n_i\}_{i=1}^\8$, $n_i \ra
\8$ as $i \ra \8$ such that $(\phi')^{n_i}(q) \ra q$ as $i\ra \8$.

\begin{proof} Let $\ep >0$ be given. For any $q\in W^u_\phi(p)$
and $U \subset M^n$, a small neighborhood of $q$, choose $\de_1$ small so
that $B_{\de_1}(q) \subset U$. Consider the small ball $B^1 =B_{c\de^r_1
  \ep_1}(q)$, where $c$ is given by the perturbation lemma and $0< \ep_1
\leq \ep /2$. Since $\phi$ preserves the volume and $M^n$ is compact, there
exists a positive integer $j_1$ such that $\phi^{j_1}(B^1) \bigcap
B_{\de_1}(q) \neq \emptyset$ and $\phi^i(B^1) \bigcap B_{\de_1}(q) =
\emptyset$ for all $i=1, 2, \ldots, j_1-1$. This implies that there exists
a point $q_1 \in B_{\de_1}(q)$ such that $\phi^{j_1}(q_1) \in B_{\de_1}(q)$
and $\phi^i(q_1) \notin B_{\de_1}(q)$, for all $i=1, 2, \ldots, j_1-1$. Now
we apply the perturbation lemma to obtain $\phi_1 \in \dfr$ with $\|\phi_1
-\phi\|_{C^r} < \ep_1 \leq \ep/2$ and Support$(\phi -\phi_1) \subset
B_{\de_1}(q)$ such that $\phi_1(q) = \phi(q_1)$.  Thus $(\phi_1)^{j_1}(q)
\in B_{\de_1}(q)$.

Now, since $(\phi_1)^{j_1}(q) \neq q$, we may choose $0 < \de_2 < \de_1/2$
so that $(\phi_1)^j(q) \notin B_{\de_2}(q)$. Let $B^2 =B_{c\de^r_2
  \ep_2}(q)$, where $0< \ep_2 \leq \ep /4$. We choose $\de_2$ and $\ep_2$
so small such that for all $\phi' \in \dfr$, with $\|\phi' -\phi\|_{C^r}
\leq \ep_2$, and any point $x \in B^2$, we have $(\phi')^{j_1} (x) \in
B_{\de_1}(q)\backslash B_{\de_2}(q)$ and $(\phi')^{i} (x) \notin
B_{\de_1}(q)$ for all $i=1,2 \ldots, j_1-1$.

With the same argument, we see that there is a point $q_2\in B^2$ and an
integer $j_2 > j_1$ such that $(\phi_1)^{j_2}(q_2) \in B_{\de_2}(q)$ and
$(\phi_1)^i(q_2) \notin B_{\de_2}(q)$, for all $i=1, 2, \ldots, j_2-1$.
Again we apply the perturbation lemma to obtain $\phi_2 \in \dfr$ with
$\|\phi_1 -\phi_2\|_{C^r} < \ep_2 \leq \ep/4$ and Support$(\phi_2 -\phi_1)
\in B_{\de_2}(q)$ such that $\phi_2(q) = \phi_1(q_2)$.  Observe that
$(\phi_2)^{j_1}(q) \in B_{\de_1}(q)$ and $(\phi_2)^{j_2}(q) \in
B_{\de_2}(q)$.

Continue the above process, we obtain a sequence of real positive numbers
$\de_1$, $\de_2$, $\ldots$, a sequence of integers $0 < j_1 < j_2 <
\ldots $, and a sequence of functions $\phi_1$, $\phi_2$, $\ldots, \in
\dfr$ such that $(\phi_i)^{j_k} \in B_{\de_k}(q)$ for all $k=1, 2, \ldots, i$.

Let $\phi' = \lim_{i\ra \8} \phi_i \in \dfr$, then $\|\phi' - \phi\|_{C^r}
< \ep$ and $q$ is a recurrent point of $\phi'$.

This proves the lemma.
\end{proof}

The above lemma can also be applied to stable branches to obtain
backward recurrent point. In fact, the proof of the above lemma yields
a stronger result, which we will state next. For every hyperbolic
periodic point on a surface, the stable manifold (and unstable
manifold) is separated by the periodic point itself into two
branches. We have the following lemma.

\begin{lem} \label{lem40}
  There is a residual subset $R \subset \dfr$ such that if $f \in R$,
  and $p$ is a hyperbolic periodic point, then every branch of stable
  and unstable manifolds of $p$ has a forward or backward recurrent
  point.
\end{lem}

As we will prove in the next sections, rational flux vector and the
existence of recurrent points on each branch of stable and unstable
manifolds imply the existence of homoclinic points.

Before we proceed, we state another lemma on accumulations of stable
and unstable manifolds.

\begin{lem}\label{lem50}
  Let $f \subset \dfr$ be an area preserving diffeomorphism and let
  $p$ be a hyperbolic periodic point. Let $B_1$ and $B_2$ be branches
  of stable or unstable manifolds of $p$. $B_1$ and $B_2$ may be the
  same branch. If $B_1 \cap L(B_2) \neq \emptyset$, then $B_1 \subset
  L(B_2)$. Where $L(B_2)$ is the $\alpha$-limit set of $B_2$ if $B_2$
  is the stable branch and $L(B_2)$ is the $\omega$-limit set of $B_2$
  if $B_2$ is the unstable branch.
\end{lem}

This follows from a theorem of Mather \cite{Mather82a}. See also
Oliveira \cite{Oliveira00} for a proof.

Let $p$ be a hyperbolic periodic point of period $k$ for $f \in \dfr$.
We may assume, by iterating $f^k$ twice, that each branch of the
stable or unstable manifold is invariant under $f^k$. We take a
linearization near $p$ such that $p$ is the origin and locally the
$x$-axis is unstable manifold, $y$-axis is the stable manifold.
We denote the stable and unstable branches by $B_u^\pm$ and
$B_s^\pm$.

We first consider the branch of unstable manifold from the positive
$x$-axis. Let $q \in B_u^+$ be a recurrent point $B_u^+$. The orbit of
$q$ accumulates on $\{q\}$ itself. This implies that the orbit of $q$
accumulates at some point on the stable manifold of $p$. The orbit of
$q$ can accumulate the $y$ axis in two ways: from the first quadrant
or from the fourth quadrant, or possibly both. We may assume that it
is the former case, the orbit of $q$ accumulates from the first
quadrant. Let $B_s^+$ be the branch of the stable manifold from the
positive $y$-axis. By Lemma \ref{lem50}, we have that $B^+_s \subset
L(B^+_u)$ and $B^+_u \subset L(B^+_u)$.

By Lemma \ref{lem40}, there is a backward recurrent point on
$B_s^+$. The orbit of this recurrent point approaches accumulate on
itself either from the first quadrant or the third quadrant. If it
is the former case, then we have an adjacent pair of stable and
unstable branches $B^+_u$ and $B^+_s$, that accumulate on each
other. See Figure \ref{hpfig1}.

We claim that there is always an adjacent pair of stable and unstable
branches accumulate on each other. In the above case, if $B_s^+$
accumulates on itself through the third quadrant, then we can consider
the unstable branch $B_u^-$ on the negative $x$-axis and the stable
branch $B_s^-$ on negative $y$-axis.  Again each branch has two ways
of accumulating on itself. Either there are two adjacent branches
accumulating on each other, or we have a cyclic accumulation: $B^+_u$
on $B^+_s$, $B^+_s$ on $B^-_u$, $B^-_u$ on $B^-_s$, $B^-_s$ and on
$B^+_u$. However, by Lemma \ref{lem50}, the above accumulations are
transitive, $B^+_s$ would have to accumulate on $B^-_s$ and hence
$B^+_u$ and therefore, $B^+_u$ and $B^+_s$ are accumulating adjacent
pair.

\section{Maps with zero flux vector}

In this section, we assume that $f \in \dfro$ is homotopic to identity
and $\bv_f =0$.

\begin{figure}[htb]
\begin{center}
\includegraphics[height=2.5in,width=2.5in]{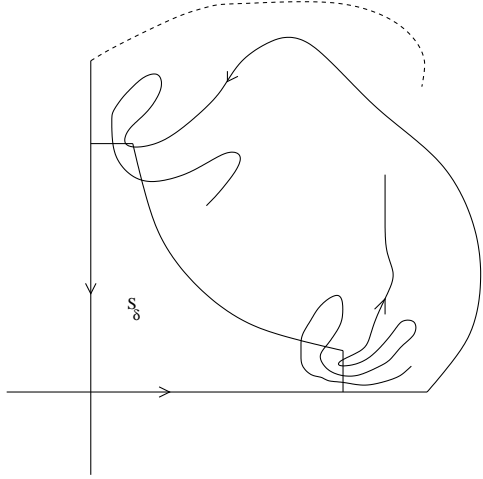}
\caption{Accumulation of stable and unstable branch} \label{hpfig1}
\end{center}
\end{figure}

Let $p$ be a hyperbolic periodic point of period $k$ and let $W^s(p)$
and $W^u(p)$ be, respectively, the stable manifold and the unstable
manifold of $p$. For simplicity and without loss of generality, we
choose a local coordinate around $p$, $$D = \{(x, y) \in \R^2 \; | \;
|x| \leq 2\lambda^k \eta, |y| \leq 2\lambda^k \eta \}$$ with some
$\lambda >0$ and some small $\eta >0$ such that $p$ is the origin in
$D$ and $f^k$ is linear $D$,
$$f^k(x, y) = (\lambda^k x, \lambda^{-k} y)$$ for some $\lambda > 1$.
Let $B_u$ be the branch of the unstable manifold of $p$ containing the
positive $x$-axis and Let $B_s$ be the branch of the stable manifold
of $p$ containing the positive $y$-axis. We further assume that $B_s$
accumulates in $B_u$ from the first quadrant in $D$ and $B_s$
accumulates in $B_u$ from the first quadrant in $D$, as
Figure~\ref{hpfig1} shows. The main result of this section is the
following theorem.

\begin{thm} \label{thm10}
  Let $f \in \dfro$ be an area preserving diffeomorphism. Assume that
  $f$ is isotopic to identity and $f$ has zero flux, $\bv_f =0$. Let
  $p$ be a hyperbolic periodic points with $B_s$ and $B_u$ as
  described above. Then $B_u \cap B_s \backslash \{p\} \neq
  \emptyset$. i.e., there are homoclinic points for the hyperbolic
  periodic point $p$.
\end{thm}

\begin{figure}[htb]
\begin{center}
\includegraphics[height=2.5in,width=3.5in]{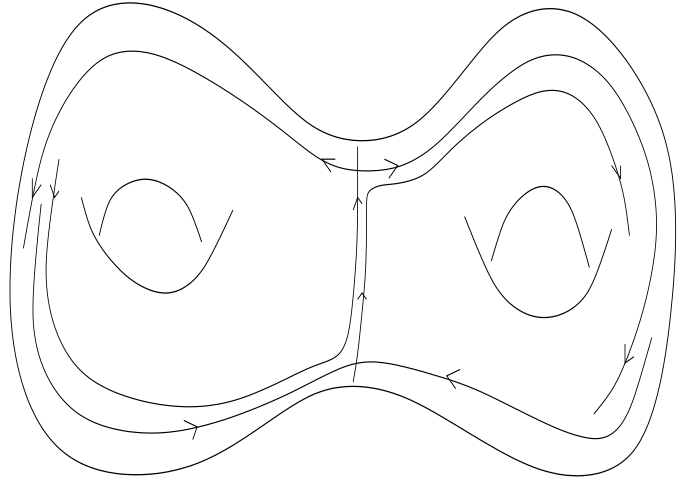}
\caption{Dense stable and unstable manifolds without homoclinic
  point.} \label{hpfig5}
\end{center}
\end{figure}

We remark that the assumption on the flux is crucial. For surface
$M_g$ with genus $g \geq 2$, one can easily construct examples of area
preserving flows with exactly $2g-2$ hyperbolic fixed points and each
branch of the stable and unstable manifolds is dense in the whole
manifold and there is neither homoclinic point nor heteroclinic
point. See Figure \ref{hpfig5}.

We will prove the theorem in a sequence of lemmas. For a positive
integer $m$, let $\delta$ be a positive real number such that
$\lambda^{km} \delta = \eta$. Let $S \subset D$ be a region defined by
$$S_\delta = \{ (x, y) \in D \; | \; xy \leq \delta \eta,\; 0\leq x \leq
\eta,\; 0 \leq y \leq \eta \}$$ 

We use the natural order for the points on $B_u$ and $B_s$. For any $z
\in B_u$, we denote the segment of $B_u$ from $p$ to $z$ by $B_u[0,
z]$. We also write $B_u = B_u[0, \8)$ and likewise, $B_s = B_s(-\8,
0]$. Let $u_1$ be the smallest point (in other words, the first point
starting from $p$) on $B_u$ that re-enters the region $S_\delta$. The
point $u_1$ is on the boundary of $S_\delta$. Moreover, it enters
either from the line $\{0 \leq x \leq \delta, \; y = \eta \}$ or from
the curve $\{ xy = \eta\delta, \; \delta \leq x \leq \lambda^k \delta
\}$. For simplicity, we assume that $u_1$ is on the line $\{ y = \eta
\}$, the other case can be dealt in a similar way.

For any two points $z_1, z_2 \in D$, let $l[z_1, z_2]$ be the oriented
straight line segment from $z_1$ to $z_2$. Obviously $B_u[0, u_1] *
l[u_1, 0]$ is an oriented simple closed curve. Similarly, we let $z_1$
be the first point on the stable branch $B_s$ that enters
$S_\delta$. Joining two line segment $B_s[0, s_1]$ and $l[s_1, 0]$ we
also obtain an oriented simple closed curve. If our manifold $M$ is
$S^2$, we already see that $B_s(0, s_1)$ and $B_u(0, u_1)$ have to
intersect, for otherwise we have two contractible simple closed curves
$B_u[0, u_1] * l[u_1, 0]$ and $B_s[0, s_1] * l[s_1, 0]$ cross each
other only once, at the origin. This is impossible.

For general surfaces, closed curves may not necessarily be
contractible. The situation is much more complicated, this simple
argument no longer work.

\begin{lem} \label{lem10}
Assume that $f \in \dfro$ is isotopic to identity and $f$ has zero
flux. Suppose that the hyperbolic periodic point has no homoclinic
point. Then there is a sequence of points $u_1 < u_2 < u_3 < \ldots,$
on $B_u$ such that

(a) $\pi_y(u_i) = \eta$ and $0 < \ldots < \pi_x(u_3) < \pi_x(u_2) <
\pi_x(u_1) \leq \delta$, where $\pi_x$ and $\pi_y$ are projections
into respective coordinates;

(b) The length of the curve $B_u[0, u_i]$, $|B_u[0, u_i]| \ra \8$ as
$i \ra \8$; and

(c) The closed curves $B_u[0, u_i]  *  l[u_i, 0]$, $i \in \N$, are
all in the same homotopy class.
\end{lem}

\begin{proof}
To simplify the notations, we assume that $p$ is a fixed point. i.e.,
$k=1$. The general case with $k >1$ is exactly the same.

Let $s=(0, \eta)$ be the point on the stable manifold $B_s$ and, as in
the statement of the lemma, $u_1$
be the first intersection of $B_u$ in $S_\delta$. The curve
$$l_1 = B_u(0, u_1)  *  l[u_1, s]  *  l[s, 0]$$
form a simple closed curve. Let $s_1 = f(s) =(0, \lambda^{-1} \eta)
\in B_s$ be the image of $s$. The the image of $l_1$, 
$$f(l_1) = B_u(0, f(u_1))  *  l[f(u_1), s_1]  *  B_s[s_1, 0]$$
is a closed curve, homotopic to $l_1$.  Moreover, since the flux of
$f$ is zero, the signed area enclosed by $l_1$ and $f(l_1)$ is
zero. This implies that the piece of the unstable manifold $B_u[u_1,
f(u_1)]$ intersects the line $l[u_1, s_1]$ at least at two distinct
points. There may be, and will be, many points of intersections $z$
such that $B[0, z]  *  l[z, 0]$ is not in the same homotopy class as
that of $B[0, u_1]  *  l[u_1, 0]$. But at least two distinct points
will have this property, since $l_1$ and $f(l_1)$ are homotopic. The
first point $u_1$ certainly has this property. Let $u_2$ be the last
point on $B[u_1, f(u_1)]$ such that $u_2$ is on the line $l(u_1, s)$
and $B[0, u_1]  *  l[u_1, 0]$ is homotopic to $B[0, u_1]  *  l[u_1,
0]$. Clearly, $u_2 \neq u_1$. A better way to understand the choice of
$u_2$ is from the universal covering space of surface. The point $u_2$
is simply the last intersection point of the proper lifts of $l_1$ and
$f(l_1)$. The lifts of $l_1$ are, of course, no longer necessarily
closed. See Figure~\ref{hpfig3}.

\begin{figure}[htb]
\begin{center}
\includegraphics[height=2.5in,width=4in]{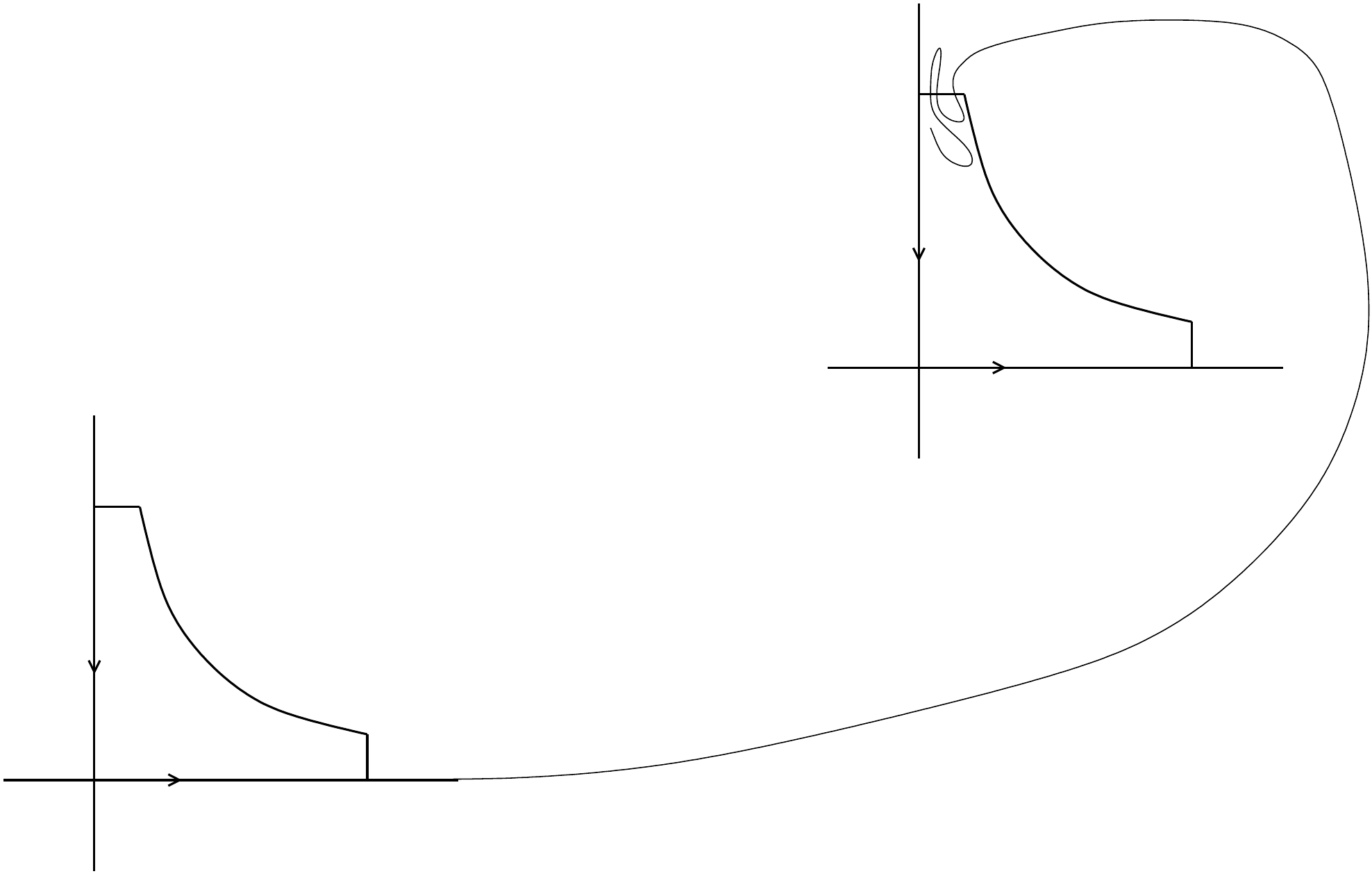}
\caption{The sequence $u_1$, $u_2$, $\cdots$} \label{hpfig3}
\end{center}
\end{figure}

We can now consider the closed curve $B_u(0, u_2)  *  l[u_2, s]  * 
B_s[s, 0]$ and its image. In the same way, we obtain an intersection
point $u_3$. Continue this process, we obtain a sequence, $u_1 < u_2 <
u_3 < \ldots$, satisfying the properties of the lemma.

This proves the lemma.
\end{proof}

Let $\{u_i\}$, $i \in \N$ be the sequence from the above lemma.  The
sequence of real numbers $\{ \pi_x(u_i) \}$, $i \in \N$ is
monotonically decreasing and bounded below by $0$, so it converges. It
is not obvious that the sequence $\{u_i\}_{i \in \N}$ should converge
to $s$ in the stable manifold. However, we will show that it always
does. We have the following lemma.

\begin{lem} \label{lem20}
  Let $u_1 < u_2 < u_3 < \ldots,$ be the sequence of points on the
  unstable manifold $B_u$ given by Lemma \ref{lem10}. Then
  $\pi_x(u_i) \ra 0$ as $i \ra \8$.
\end{lem}

\begin{figure}[htb]
\begin{center}
\includegraphics[height=2.0in,width=2.5in]{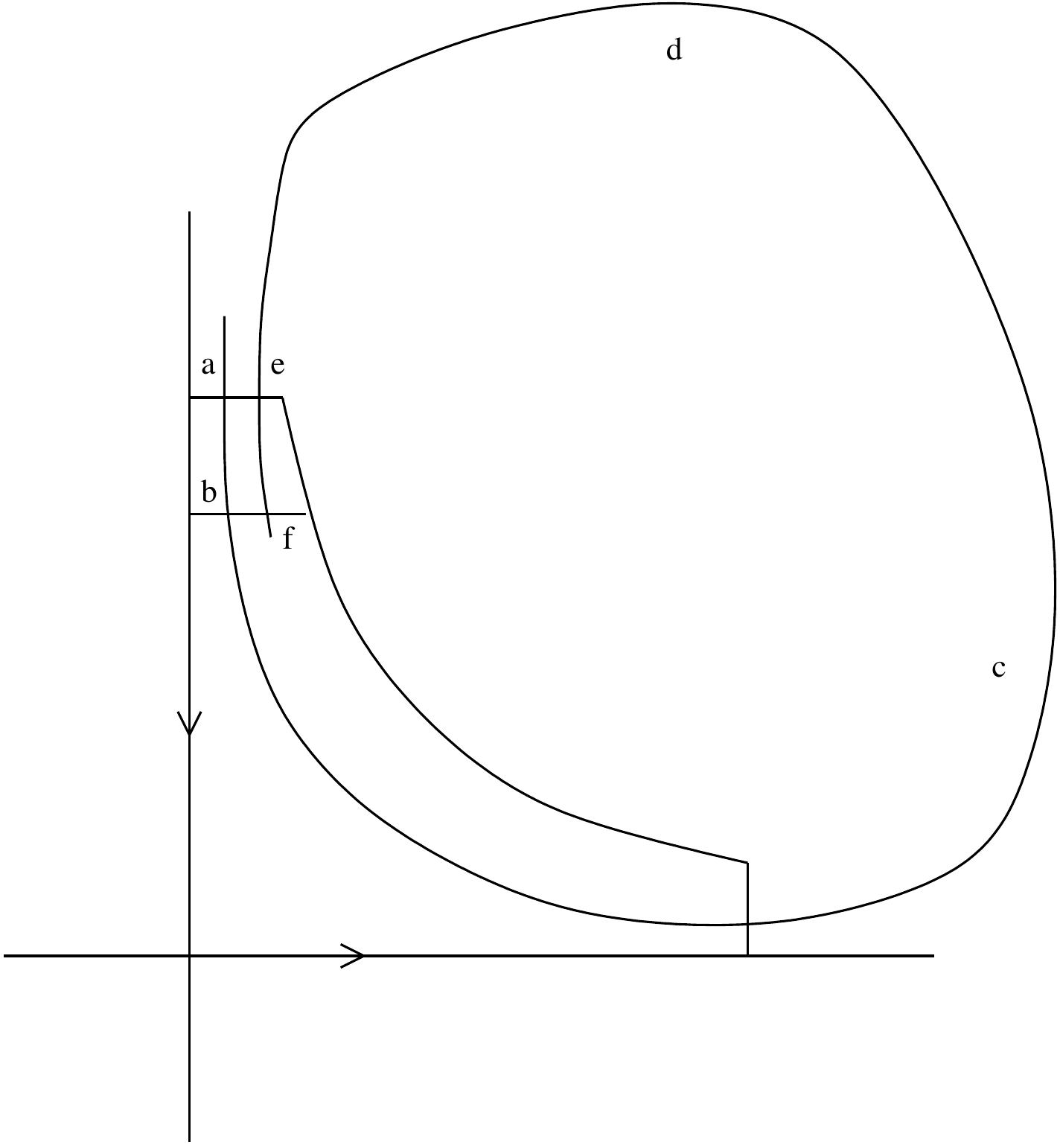}
\caption{Nonzero flux} \label{hpfig2}
\end{center}
\end{figure}

The proof of this lemma is quite technical, we postpone it to the next
section. Here we explain some basic ideas involved in the
proof. Suppose that there is curve $l$ such that it enters the region
$S_\delta$ twice at different points, as shown in
Figure~\ref{hpfig2}. Further suppose that the curve is invariant in
the sense that the segment $(abcde)$ is mapped to $(bcdef)$. We can
join two points on the curve by a horizontal line to form a closed
curve $(abcdefa)$ as in the figure. Since the curve is invariant, this
closed curve has nonzero flux, equal to the area of $(abfea)$, which
is impossible for our map. Suppose the sequence $u_i$, $i \in \N$ does
not converge to the stable manifold, then the limit is in an invariant
set formed by taking the limit of pieces of unstable manifold. Since
the unstable manifold accumulates on the stable manifold, this
invariant set can get arbitrarily close to the stable manifold. If
this invariant set is simply a curve, then above argument can derive a
contradiction. However, the limit set of a sequence of curves can be
very complicated, we can not directly apply the above argument. But
nevertheless, we can only use curves to approximate the invariant set,
as long as we can control the area lost by the approximation.

\vs{1ex}
\noindent {\bf Proof of Theorem \ref{thm10}}: We are ready to prove
the main theorem of this section. 

By Lemma \ref{lem20}, there is a sequence of points $u_i$, $i \in \N$
on the unstable manifold $B_u$ such that $u_i \ra (0, \eta)$ and
$B_u[0, u_i] * l[u_i, 0]$ is in the same homotopy class for all $i \in
\N$. Likewise there is a sequence $s_i$, $i \in \N$ in the stable
manifold $B_s$ such that $s_i \ra (\eta, 0)$ and $B_s[s_i, 0] * l[0,
s_i]$ are in the same homotopy class for all $i \in \N$. For any $\ep
>0$, let $$S_\ep = \{ (x, y) \in D \; | \; xy \leq \ep \eta,\; 0\leq x
\leq \eta,\; 0 \leq y \leq \eta \}$$ We claim that there is point $z$
in $B_u$ such that $z \in S_\ep$ and the closed curve $B_u[0, z] *
l[z, 0]$ is homotopic to the closed curve $$B_u[0, u_1] * l[u_1, 0] *
B_s[0, s_1] * l[s_1, 0].$$ This can be shown by the following: Let
$z'$ be a point in $B_s$ which is close to some $s_i$ and inside of
$S_\ep$. $B_s[0, z'] * l[z'. 0]$ is homotopic to $B_s[0, s_i] *
[s_i, 0]$. There exists an integer $j$ such that $f^{j}(z')$ is
between the points $(0, \eta)$ and $(0, \lambda^{-1}\eta)$ on the
stable branch $B_s$. By choosing a different $\eta$, we may just
assume that $f^j(z') = (0, \eta)$. For large $i$, the point$u_i$ is
close to $(0, \eta)$, therefore we have that $f^j(u_i) \in
S_\ep$. Since $B_u[0, u_i] * l[u_i, (0, \eta)] * l[(0, \eta), 0]$ is
homotopic to $B_u[0, u_1] * l[u_1, 0]$, so is $f^{-j}(B_u[0, u_i] *
l[u_i, (0, \eta)] * l[(0, \eta), 0])$. This implies that $B_u[0,
f^{-j}(u_i)] * l[f^{-j}(u_i), 0]$ is homotopic to $B_u[0, u_1] *
l[u_1, 0] * B_s[0, s_1] * l[s_1, 0]$. We may choose $z =
f^{-j}(u_i)$. This proves the claim.

In fact, for later convenience, we will choose the point $z$
differently. The trajectory of $f^{-j}(u_i)$ enters the set $S_\ep$
before $f^{-j}(u_i)$. We let $z$ be the first point on $B_u$ such that
$z$ is in the set $S_\ep$ and $B_z[0, z]*l[z, 0]$ is homotopic to
$B_z[0, f^{-j}(u_i)]*l[f^{-j}(u_i), 0]$. Obviously, $z$ must be either
on the line $\{0 < x < \ep, \; y = \eta \}$ or on the curve $\{xy =
\ep \eta, \; \lambda^{-1} \eta \leq y \leq \eta \}$.

Similarly, there is a point $w \in B_s$ such that $w$ is the first
point in $B_s$ to intersect $S_\ep$ such that the closed curve $B_s[0,
w] * l[w, 0]$ is homotopic to the closed curve $B_s[0, s_1] * l[s_1,
0] * B_u[0, u_1] * l[u_1, 0]$. Similarly, $w$ must be either
on the line $\{0 < y < \ep, \; x = \eta \}$ or on the curve $\{xy =
\ep \eta, \; \lambda^{-1} \eta \leq x \leq \eta \}$.

We have obtained two closed curves $B_u[0, z] * l[z, 0]$ and $B_s[0,
w] * l[w, 0]$. They are not homotopic to each other, but they are
homologous to each other, since $B_u[0, u_1] * l[u_1, 0] * B_s[0, s_1]
* l[s_1, 0]$ and $B_s[0, s_1] * l[s_1, 0] * B_u[0, u_1] * l[u_1, 0]$
are homologous. These two curves cross each other at the origin. This
crossing at the origin has an intersection number $\pm 1$. The sign in
$\pm 1$ depends on the orientation we pick for these two curves. There
are possibly many other intersection between $B_s[0, w] * l[w, 0]$ and
$B_u[0, z] * l[z, 0]$. However, we want to show that all other
intersections contribute to a total of zero intersection number. Since
any tow homologous curves have a total of zero intersection number, we
reach a contradiction and hence the theorem is proved.

Let $z_1 =z$, apply Lemma \ref{lem10} and Lemma \ref{lem20} to the
point $z_1$, we obtain a sequence of points $z_i$, $i \in \N$ on the
line $y = \eta$ such that the following is true:
\begin{enumerate}
\item the sequence $z_i$, $i \in \N$ is on the unstable branch $B_u$
  and $z_1 < z_2 < z_3 < \cdots$ with the order on the unstable
  branch.

\item $\pi_x(z_1) > \pi_x(z_2) > \cdots$ and $\lim_{i \ra \8}
  \pi_x(z_i) = 0$.

\item For any $i \in \N$, the closed curve $B_u[0, z_i]*l[z_i, 0]$ is
  homotopic to $B_u[0, z]*l[z, 0]$.
\end{enumerate}
Similarly, we obtain a sequence of points $\{w_i\}_{i \in \N}$ on the
stable branch $B_s$ with $w_1=w$.

For any $i_0 \in \N$, the curve $B_u[0, z_{i_o}]$ is finite and
therefore does not intersect with $l[w_i, 0]$ for sufficiently large
$i$. For small $i$, the curve $B_u[0, z_{i_o}]$ may have intersection
with $l[w_i, 0]$.  However, the curve $l[0, w_i]*B_s[w_i, w_j]*l[w_j,
0]$ is a simply closed curve, therefore the total intersection number
between $B_u[0, z_{i_o}]$ and $l[w_i, 0]$ is zero for any $i$. In
other words, the total intersection number between $B_u[0, z]$ and
$l[w, 0]$ is zero. Likewise the total intersection number between
$B_s[0, w]$ and $l[z, 0]$. This shows that the total intersection
number between $B_u[0, z] * l[z, 0]$ and $B_s[0, w] * l[w, 0]$ is
$\pm1$, contradicting to the fact that these to curves are
homologous. This contradiction show that the stable and unstable
branches $B_s$ and $B_u$ have to intersect.

This proves the theorem.

\section{Proof of Lemma \ref{lem20}}

In this section, we give a proof of Lemma \ref{lem20}.  We will prove
by contradiction. Suppose on the contrary that $\lim_{i \ra \8}
\pi_x(u_i) = q_x >0$. We write $q =(q_x, \eta) = \lim_{i \ra \8} u_i$.

Let $\ep$ be a positive real number such that $\ep < q_x/2$ and let
$$S_\ep = \{ (x, y) \in D \; | \; xy \leq \ep \eta,\; 0\leq x \leq
\eta,\; 0 \leq y \leq \eta \}$$ Since the unstable branch $B_u$
accumulates on the stable branch $B_s$ from the first quadrant, there
is a point $v_1$ on $B_u$ such that $v_1$ is the first intersection of
$B[0, v_1]$ with $S_\ep$. Again, we assume that $v_1$ is on the line
$\{0 \leq x \leq \delta, \; y = \eta \}$. The case where $v_1$ is on
the curve $\{ xy = \eta\ep, \; \ep \leq x \leq \lambda \ep \}$
can be dealt in a similar way.

Similar to the choices of the sequence $\{u_i\}_{i \in \N}$ in Lemma
\ref{lem10}, we obtain a sequence of $\{v_i \}_{i \ra \8}$ having the
same properties as in Lemma \ref{lem10}, except of course that $B_u[0, v_i]
 *  l[v_i, 0]$ is in a different homotopy class as that of $B_u[0, u_i]
 *  l[u_i, 0]$.

For each $v_i$, $i \in \N$, there is a point $u_i' \in B_u$, $u_i' \in
B_u[0, v_i]$ such that
\begin{enumerate}
\item $\pi_y(u_i') = \eta$ and $0 < \pi_x(u_i') < \pi_x(u_1)$;

\item $B_u[0, u_i']  *  l[u_i', 0]$ is in the same
homotopy class as $B_u[0, u_i]  *  l[u_i, 0]$;

\item  the point $u_i'$ is the last such point on $B_u[0, v_i]$ to
  have the above properties. 
\end{enumerate}

\begin{figure}[htb]
\begin{center}
\includegraphics[height=2.5in,width=3.5in]{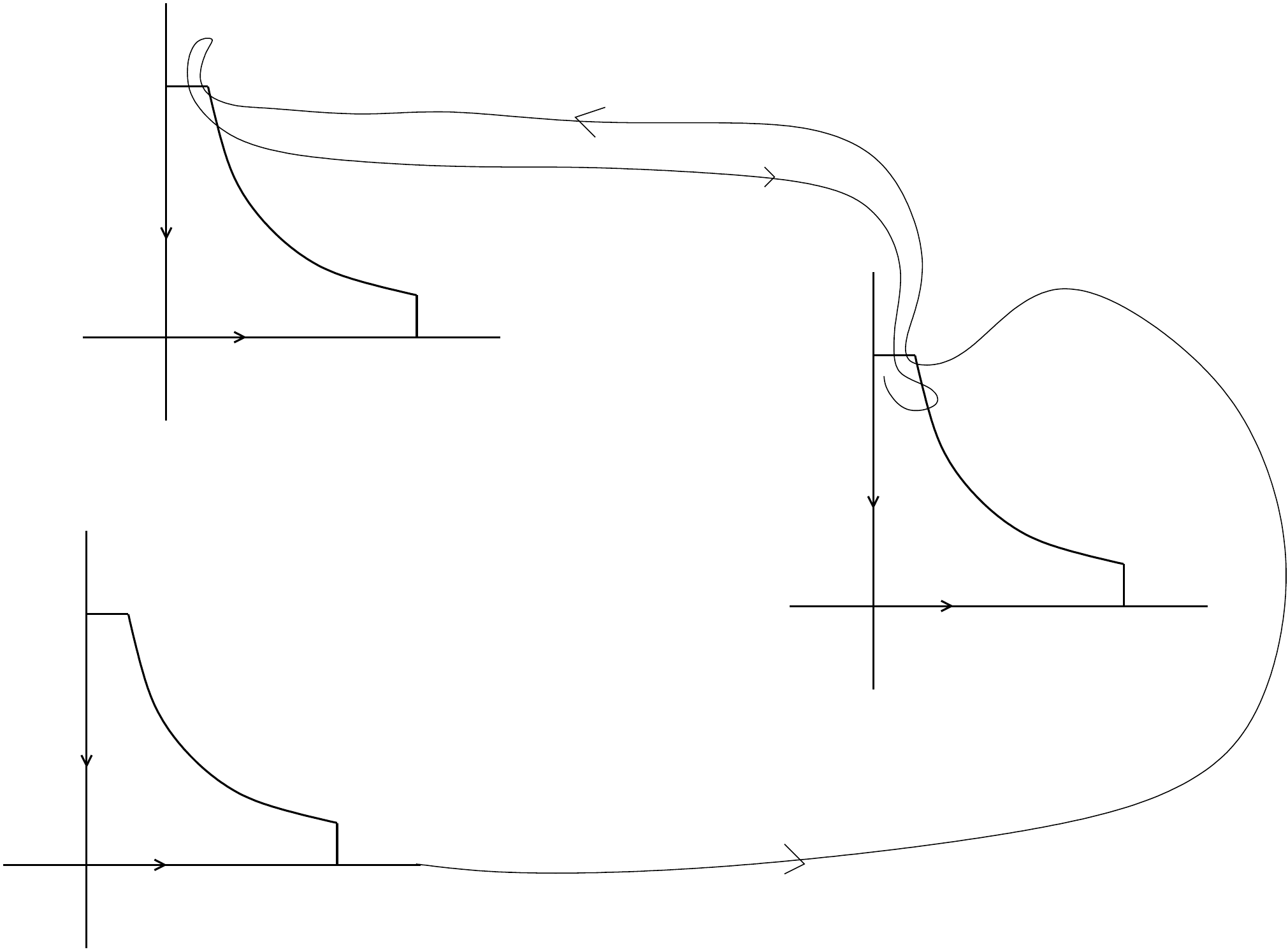}
\caption{Sequences $u_i'$ and $v_i'$.} \label{hpfig4}
\end{center}
\end{figure}

Again, the choice of $u_i'$ is easy if one looks from the universal
covering space. See Figure~\ref{hpfig4}. Finally we choose $v_i'$ to
be the first point in $B_u[u_i', v_i]$ and on the line $\{y = \eta, 0
< x < \ep \}$ such that $B_u[0, v_i'] * l[v_i', 0]$ is homotopic to
$B_u[0, v_i] * l[v_i, 0]$.

Now we have two sequences of points, $\{u_i'\}_{i \in \N}$ and
$\{v_i'\}_{i \in \N}$. To simplify our notations, we will assume that
our original choices of $u_i$ and $v_i$ already satisfy the properties
for $u_i'$ and $v_i'$ we have stated.

We let $q' = (q_x', \eta) = \lim_{i \ra \8} v_i$. We have $0 \leq q_x'
< \pi_x(v_i)$, for $i \in \N$. 

Now we have a sequence of lines $B_u[u_i, v_i]$, $i \in \N$. Every two
consecutive lines $B_u[u_i, v_i]$ and $B_u[u_{i+1}, v_{i+1}]$,
together with the line segments $l[u_i, u_{i+1}]$ and $l[v_i,
v_{i+1}]$ bound an area. On the universal cover with a fixed lift, the
area bounded by these four curves do not have self intersections and
it doesn't intersect other such areas, except at the
boundary. However, when these areas are projected down to the
manifold, there might be multiple covers at some open sets. Since we
will use the finiteness of the total area to show certain convergence
properties of the sequence of the curves $B_u[u_i, v_i]$, We need to
changes the curves $B_u[u_i, u_{i +1}]$ so that there will be no
multiple covers. We will do it in the following way.

For any $i > 1$, let $A_i$ be the strip bounded by $B_u[u_i, v_i]$,
$B_u[u_{i+1}, v_{i+1}]$, $l[u_i, u_{i+1}]$ and $l[v_i, v_{i+1}]$ and
let $A = \cup_{i \in \N} A_i$. We will construct a sequence of new
curves $d_i$, $i \in \N$, to replace $B_u[u_i, v_i]$. In the universal
covering space, the strip $A_i$ and the set $A$ have infinitely many
lifts. Fix one of the lift of $A_i$, say, $\tilde{A}_i$. The set
$\tilde{A}_i$ may intersect other lifts of itself, or lifts of other
strips, $A_j$, $j \in \N$ in the interior. If this does not happen, we
will simply let $d_i = B_u[u_i, v_i]$ and $d_{i+1} = B_u[u_{i+1},
v_{i+1}]$. Otherwise, we will remove the interior of any part of
$\tilde{A}_i$ that intersects with any other lift of $A_i$ and any
lifts of $A_j$, $j \neq i$.

To understand the structure of the sets removed from $A_i$, the
following fact is useful: if any curve $B_u[u_i, v_i]$ enters a strip
$A_j$, $j \neq i$, it has to enter through the end of the strips,
either through the line segment $l[u_j, u_{j+1}]$ or the line segment
$l[v_j, v_{j+1}]$. Moreover, it has to eventually exit through the
same segment. This is from our construction that $B_u[u_i, v_i]$
contains no proper segment crossing the strips.

On the manifold, the areas we removed therefore are disks with pieces
of boundaries on the line $\{y = \eta, 0 < x < u_1 \}$. We now let
$d_i$ and $d_{i+1}$ be the new boundary pieces of $A_i$. The new
boundaries $d_i$ and $d_{i+1}$ are homotopic to, respectively,
$B_u[u_i, v_i]$ and $B_u[u_{i+1}, v_{i+1}]$, relative to their end
points.

We now have a sequence of curves $d_i$, $i \in \N$ and a sequence of
strips $B_i$ on the surface, bounded by $d_i$, $d_{i+1}$, $l[u_i,
u_{i+1}]$ and $l[v_i, v_{i+1}]$. The sequence of strips $B_i$ does not
intersect each other except at the boundaries. Since the total area of
the surface is bounded, we conclude that the area of $B_i$ approaches
zero, as $i \ra \8$. This will be an important fact.

We now consider the image of $d_i$ under the map $f$. Recall that by
our assumption, $\lim_{i \ra \8} u_i = q=(q_x, \eta)$ and $\lim_{i \ra
  \8} v_i = q=(q_x', \eta)$. Reduce $\delta$ if necessary, we observe
that the curves $d_i$ has to either intersect the line $y = \lambda
\eta$ or the line $y = \lambda^{-1} \eta$, right after the point
$u_i$. For definiteness, we may, by choosing a subsequence if
necessary, assume that all $d_i$ intersect the line $y = \lambda^{-1}
\eta$ right after the point $u_i$. With this choice, by the
orientation preserving property, $d_i$ intersects the line $y =
\lambda \eta$ right before the point $v_i$.

Let $u_i^1$ be the intersection point of $d_i$ with the line $y =
\lambda^{-1} \eta$ right after $u_i$. i.e., $u_i^1$ is the first point on
$B_u[u_i, v_i]$ that is on $y = \lambda^{-1} \eta$. Clearly, $f(u_i)$
is on the right hand side of $u_i^1$, but left hand side of
$u_{i-1}^1$. This implies that $$\lim_{i \ra \8} u_i^1 = \lim_{i \ra
  \8} f(u_i) = f(q).$$

This suggests certain invariance properties of the limit set of the
curves $B_u[u_i, v_i]$ and $d_i$. Indeed, let $z$ be a limit point of
the curves $B_u[u_i, v_i]$. i.e., there is a sequence of the points
$z_i \in B_u[u_i, v_i]$ such that $\lim_{i \ra \8} z_i =z$, then
$f(z_i)$ is between two curves $B_u[u_i, v_i]$ and $B_u[u_{i-1},
v_{i-1}]$, as long as $f(z_i)$ stays inside of the strip $A_i$ bounded
by $B_u[u_i, v_i]$, $B_u[u_{i+1}, v_{i+1}]$, $l[u_i, u_{i+1}]$ and
$l[v_i, v_{i+1}]$. Let $z_i'$ be a point on the intersection of the
line $l[f(z), f(z_i)]$ and the curve $B_u[u_i, v_i]$, we have $f(z) =
\lim_{i \ra \8} z_i'$., i.e., $f(z)$ is also a limit point of the
curves $B_u[u_i, v_i]$.

This shows that the limit set of the sequence of curves $B_u[u_i,
v_i]$ is forward and backward invariant except near the ends of the
curves. This limit set is connected and, in general, it is expected to
be extremely complicated. In certain sense, this limit set starts at
$q$ and ends at $q'$, the limit of the sequence $v_i$. Near these two
end points, the limit set is quite easy to describe. Between the two
points $q$ and $f(q)$, this limit set is simply the limit of the
curves $B_u[u_i, u_i^1]$. Similarly, between $f^{-1}(q')$ and $q'$,
this limit set is the limit of the curves $B_u[v_i^1, v_i]$, where
$v_i^1$ is defined in the similar way as $u_i^1$.

For any $\ep' >0$, there is a integer $i_0 >0$ such that, if $i \geq
i_0$, the area of the strip $B_i$ bounded by $d_i$, $d_{i+1}$, $l[u_i,
u_{i+1}]$ and $l[v_i, v_{i+1}]$ is smaller than $\ep'$. Let $E_i$ be
the square bounded by $B_u[u_i, u_i^1]$, $B_u[u_{i+1}, u_{i+1}^1]$,
$l[u_i, u_{i+1}]$ and $l[u_i^1, u_{i+1}^1]$ and let $E = \cup_{i=1}^\8
E_i$. Similarly, let $E_i'$ be the square bounded by $B_u[v_i^1,
v_i]$, $B_u[v_{i+1}^1, v_{i+1}]$, $l[u_i, u_{i+1}]$ and $l[u_i^1,
u_{i+1}^1]$ and let $E' = \cup_{i=1}^\8 E_i'$. By dropping first few
terms of the sequence $u_i$ and $v_i$, we may assume that the area of
$E$ and $E'$ are both smaller than $\ep'$.

Next, we consider the image of $B_u[u_{i+1}, v_{i+1}]$ under $f$ from
the universal covering space. For most part, the image is in the strip
$A_i$, except near the end point $v_{i+1}$, where it goes out of the
strip $A_i$. The curve $d_{i+1}$ is the shortened version of
$B_u[u_{i+1}, v_{i+1}]$, its image are mostly in the strip $B_i$,
bounded by $d_i$, $d_{i+1}$, $l[u_i, u_{i+1}]$ and $l[v_i, v_{i+1}]$,
except near the cut off points, where the image goes out of the strip
$B_i$ near $v_j$, $j \in \N$.

Now, we consider the close curve $C_{i+1} = d_{i+1 } *  l[u_{i+1},
v_{i+1}]$. The image of the curve, $f(C_{i+1})$ is homologous to
$C_{i+1}$, together they bound an oriented disk chain, say
$D_{i+1}$. i.e., $\partial D_{i+1} = f(C_{i+1}) - C_{i+1}$. By the zero
flux property, this disk chain has total area zero. We will show that
this is impossible, hence a contradiction.

The disk chain $D_{i+1}$ can be divided into several parts. The first
part is inside the strip $B_i$. Let $D_{i+1}^1 = D_{i+1} \cap B_i$,
this covers the part of $f(C_{i+1})$ inside the strip $B_i$. For
choices of large $i$, we have that
$$|\int_{D_{i+1}^1} d\mu | < | \int_{B_i} d\mu | < \ep'.$$
The second part of $D_{i+1}$ covers the part of $f(C_{i+1})$ near the
cut off points of $d_{i+1}$. Let
$D_{i+1}^2 = f(E') \cap D_{i+1}$, where $E'$ is the union of the
squares near $v_i$. We also have that 
$$|\int_{D_{i+1}^2} d\mu | < | \int_{E'} d\mu | < \ep'.$$
Let $D_{i+1}^3$ be the rest of $D_{i+1}$, i.e., $D_{i+1}^3 = D_{i+1}
\backslash (D_{i+1}^1 \cup D_{i+1}^2)$. The set $D_{i+1}^3$, plus or
minus an area of size $\ep'$, contains the square bounded by
$B_u[u_{i+1}, u_{i+1}^1]$, $B_u[f(v_{i+1}^1, v_{i+1}]$, $l[u_{i+1},
f(v_{i+1}^1)]$ and $l[u_{i+1}^1, f(v_{i+1})]$. Therefore, there is a
constant $m >0$, independent of $\ep'$ such that if $i$ is large
enough, $$\int_{D_{i+1}^3} d\mu  > m >0$$ This implies that the
absolute value of signed area of $D_{i+1}$ is larger than $m - 2 \ep'
>0$, if $\ep'$ is chosen sufficiently small. This contradicts to the
fact that every closed curve has zero flux.

This contradiction proves the lemma.
\qed

\section{Hamiltonian flows and the proof of the main theorem}

First, we give a proof of Theorem \ref{mainthm}.

Let $\dfro$ be the set of area-preserving diffeomorphism on $M$ that
are homotopic to identity and let $f \in \dfro$. Let $p$ be a
hyperbolic periodic point. By Lemma \ref{lem70}, there is $f_1 \in
\dfro$, $C^r$ close to $f$ such that $p_1$, slighted perturbed from
$p$, is a periodic point for $f_1$ with the same period and the flux
vector for $f_1$ is rational. By Lemma \ref{lem2}, there is map $f_3$,
$C^r$ close to $f_2$, such that each branch of stable manifold and
unstable manifold of $p_2$ has a forward or backward recurrent
point. This implies that for any positive integer $k >0$, the stable
and unstable branches accumulate on stable and unstable branches under
$f_3^k$, even though the original recurrent point may no longer be
recurrent. Moreover, there is an adjacent pair of stable and unstable
branches accumulating on each other. Let $k$ be the positive integer
such that the flux for $f_3^k$ is zero. By Theorem \ref{thm10}, the
stable and unstable manifolds of $p_2$ under $f_3^k$ intersects and
this intersection can be made transversal by an arbitrary $C^r$ small
perturbation to $f_3$. This implies that there is an open set of
diffeomorphisms in $\dfro$, arbitrarily $C^r$ close to $f$, such that
the perturbed periodic point from $p$ has a transversal homoclinic
point. Since $C^r$ generically in $\dfro$ there are only countably many
periodic point, this implies that there is a residual set $R_1 \in
\dfro$ such that for every $\phi \in \R_1$ and every hyperbolic
periodic point for $\phi$, there is a homoclinic point.

This proves Theorem \ref{mainthm} \qed

To prove Theorem \ref{mainthmh}, we need some preliminary results on
Hamiltonian flow. Let $H: M \times S^1 \ra \R$ be a $C^{r+1}$ time
periodic Hamiltonian function on $M$. Let $\phi: M \times \R \ra M$,
be the Hamiltonian flow on $M$, given by the $C^r$ Hamiltonian vector
field $JdH$. We write $\phi_t(*) = \phi(*, t)$. For any $t_0 \in S^1$,
let $f = \phi_1: M \times \{t_0\} \ra M \times \{t_0\}$ be the
Poincar\'e map of the Hamiltonian flow. We will identify $M \times
\{t_0\}$ with $M$ in the obvious way and regard $f$ as a map on
$M$. Such map $f$ is called {\it Hamiltonian diffeomorphism}\/.
Obviously, Hamiltonian diffeomorphisms preserve symplectic forms and
therefore, on surfaces, they are area-preserving. In fact we can say
more: All Hamiltonian diffeomorphisms have zero flux.

\begin{lem}
Let $f$ be a Hamiltonian diffeomorphism. i.e., $f$ is the time-one map
of a 1-periodic Hamiltonian flow on compact surface $M$, then the flux
vector for $f$ is zero.
\end{lem}

\begin{proof}
Let $H$ be a 1-periodic Hamiltonian function on $M$ and let $\phi_t$ be
its Hamiltonian flow such that $f = \phi_1$.
Let $l$ be a simple closed curve on $M$. For any real number $s >0$, let
the set $D \subset M \times \R$ be a flow tube defined by
$$D = \{ (x, t) \in M \times \R \; \| \; x = \phi_t(x_0),  0 \leq t
\leq s \}$$ By the definition of Hamiltonian flow, we have, $$\int_D
\omega -dH \wedge dt =0$$ Where $\omega$ is the symplectic form on
$M$. Let $\pi_M(D)$ be the natural projection of $D$ into $M$. For $s$
small, this projection is in a small neighborhood of the curve $l$,
there is, locally, a one-form $\alpha$ such that $d\alpha =
\omega$. We remark that this one-form is only defined locally. Such
one-form does not exist globally. By Stokes' theorem, we have \be 0
&=& \int_D \omega -dH \wedge dt = \int_{\phi_s(l) \times \{s \} - l
  \times \{0 \}} \alpha -H \wedge dt \nn \\
&=& \int_{\phi_s(l) \times \{s \} - l \times \{0 \}} \alpha =
\int_{\phi_s(l) - l} \alpha = \int_{\pi_M(D)} \omega \nn \ee 
This is true for any simple closed curve $l$. One can reach $s=1$ in
finitely many steps. This shows that the flux across $l$ is zero.

This proves the lemma.
\end{proof}

The next question is whether we can perturb the Hamiltonians to create
recurrent points on the stable and unstable manifolds. This is the same
as asking whether the perturbation lemma (Lemma \ref{lem3}) is valid
for Hamiltonian diffeomorphisms. This is an easy and standard exercise
using an flow box along a segment in a giving trajectory. We omit the
details.

The proof of Theorem \ref{mainthmh} now follows from that of Theorem
\ref{mainthm}.

%\bibliographystyle{plain}

%\bibliography{mybib}

\end{document}